\documentclass{amsart}
\usepackage{righttag}
\usepackage{epsf}
\def\hepsffile{\leavevmode\epsffile}

\theoremstyle{plain}
\newtheorem{thm}{Theorem}[subsection]
\newtheorem{cor}[thm]{Corollary}
\newtheorem{lem}[thm]{Lemma}

\theoremstyle{remark}

\theoremstyle{definition}

\newtheorem{emf}[thm]{}

\def\id{\protect\operatorname{id}}

\def\St{\protect\operatorname{St}}

\def\pr{\protect\operatorname{pr}}

\def\Z{{\mathbb Z}}

\def\R{{\mathbb R}}

\def\1{\hbox{\rm\rlap {1}\hskip.03in{\textrm I}}}
\def\Bbbone{{\rm1\mathchoice{\kern-0.25em}{\kern-0.25em}
	{\kern-0.2em}{\kern-0.2em}I}}

\def\p{\partial}

\begin{document} 
\centerline{{\em This paper will appear in Math. Scand.,\/}} 
\centerline{{\em probably in Vol {\bf 86}
(no. 1) (2000), to be issued in April 2000\/}}

\title[Homotopy groups of the space of curves on a surface]
{Homotopy groups of the space of curves on a surface.}

\author[V.~Tchernov]{Vladimir Tchernov}
\address{DMATH G-66.4, 
Eidgen\"osische Technische Hochschulle, CH-8092 Z\"urich, Switzerland} 
\email{chernov@math.ethz.ch}


\begin{abstract}
We explicitly calculate the fundamental group of the space $\mathcal F$ of all
immersed closed curves on a surface $F$. It is shown that  
$\pi_n(\mathcal F)=0$, $n\geq
2$, for $F\neq S^2, \R P^2$. It is also proved that 
$\pi_2(\mathcal F)=\Z$, and 
$\pi_n(\mathcal F)=\pi_n(S^2)\oplus\pi_{n+1}(S^2)$, $n\geq 3$, 
for $F$ equal to $S^2$ or $\R P^2$.
\end{abstract}

\maketitle 


By a surface we mean any smooth two-dimensional manifold.

\section{Introduction}
Recently the space of closed curves on
a surface attracted a lot of attention. 
The interest was initiated by the work of
V.~Arnold~\cite{Arnold}, who
axiomatically defined invariants  $\St$ and $J^{\pm}$  of generic curves on
$\R^2$. 

In order to define axiomatically this kind of invariants on an
arbitrary surface $F$ one has to know the fundamental group of the space $\mathcal
F$ of all immersed closed curves on $F$. 
However, as far as I know, this group is not calculated in the literature.
In this paper we explicitly calculate it. 
The knowledge of its properties allowed me~\cite{Tchernov} 
to generalize in a natural way
Arnold's invariants to the case of generic curves on 
an arbitrary surface. (The results of the well known paper by
S.~Smale~\cite{Smale}, where he calculated homotopy groups of the space of 
all immersed 
closed curves with the fixed initial point and the velocity vector at it, are
not sufficient for this purpose.) 

When the work described in this paper was complete and submitted to 
Mathematica Scandinavica, I received a preprint of
A.~Inshakov~\cite{Inshakov} containing similar results obtained by him
independently. (Later the preprint of Inshakov was broken into two
parts~\cite{Inshakovp1} and~\cite{Inshakovp2}.)

\centerline{\em Acknowledgments}
I am deeply grateful to Oleg Viro and Tobias Ekholm for 
many enlightening discussions.

\section{Main results}
\subsection{Basic definitions.}\label{basic}

A {\em curve\/} is a smooth immersion of an oriented circle $S^1$ into a
(smooth) surface $F$. For a surface $F$ we denote by $\mathcal F$ the space of
all curves on $F$.

Two curves $s_0$ and $s_1$ are said to be {\em regularly homotopic\/}, if
there exists a homotopy $H:S^1\times I\rightarrow F$ such that 
$H(t\times 0)=s_0(t)$, $H(t\times 1)=s_1(t)$, and $H(\bullet \times x)$ is an
immersion for every $x\in I$. This means that $s_0$ and $s_1$ are in the
same connected component of $\mathcal F$.

Two (oriented) curves with a tangency point, at which the velocity 
vectors of the two curves are 
pointing in the same direction, are said to be {\em direct tangent\/} 
to each other at this point.  

For a surface $F$ we denote by $STF$ the spherical tangent bundle of $F$
and by $\pr:STF\rightarrow F$ the corresponding locally trivial
$S^1$-fibration.

For a curve $\xi$ on $F$ we denote by $\vec \xi$ its lifting to $STF$, which
maps every point $t\in S^1$ to the direction of the velocity vector of $\xi$ at
$t$. 

We fix a point $a$ on $S^1$.
Then a curve $\xi$ represents an element of $\pi_1(F, \xi(a))$, and  
$\vec \xi$ represents an element of $\pi_1(STF, \vec\xi(a))$. 
When there is no ambiguity, we denote these two elements by $\xi$ 
and $\vec\xi$ respectively. 

\subsection{Fundamental group of the space of curves on an orientable 
surface.}
For orientable surfaces the group $\pi_1(\mathcal F,
\xi)$ appears to be much simpler than for nonorientable surfaces.

\begin{thm}\label{pi1S2}
Let $F=S^2$ and let $\xi$ be a curve on $S^2$. Then $\pi_1(\mathcal F,
\xi)=\Z_2$.
\end{thm}

\begin{thm}\label{pi1T2}
Let $F=T^2$ (torus) and let $\xi$ be a curve on $T^2$. Then $\pi_1(\mathcal F,
\xi)=\Z\oplus\Z\oplus\Z$.
\end{thm}

\begin{thm}\label{pi1orient}
Let $F\neq S^2, T^2$ be an orientable surface (not necessarily compact) 
and let $\xi$ be a curve on $F$.

\textrm{I.} If $\xi$ represents a homotopically nontrivial loop on $F$, then 
$\pi_1(\mathcal F, \xi)=\Z\oplus\Z$.

\textrm{II.} If $\xi$ represents a homotopically trivial 
loop on $F$, then $\pi_1(\mathcal
F, \xi)=\pi_1(STF)$.

\end{thm}
The proofs of Theorems~\ref{pi1S2},~\ref{pi1T2} and~\ref{pi1orient} are given
in Subsections~\ref{pfpi1S2},~\ref{pfpi1T2}, and~\ref{pi1arbitrary}, respectively.

\subsection{Fundamental group of the space of curves on a nonorientable
surface}

\begin{thm}\label{pi1RP2}
Let $F=\R P^2$ and let $\xi$ be a curve on $\R P^2$. Then $\pi_1(\mathcal F,
\xi)=\Z_4$.
\end{thm}

\begin{thm}\label{pi1K}
Let $F=K$ (Klein bottle) and let $\xi$ be a curve on $K$. 

\textrm{I.} If $\xi$ represents an orientation preserving loop on $K$, then 
$\pi_1(\mathcal F, \xi)=\pi_1(STK)$, provided that 
$\vec \xi=b^{2l}$ in $\pi_1(STK, \vec\xi(a))$ 
for some $b\in\pi_1(STK, \vec\xi(a))$ projecting to an orientation
reversing loop on $K$, and 
$\pi_1(\mathcal F, \xi)=\Z\oplus\Z\oplus\Z$ otherwise. 

\textrm{II}. If $\xi$ represents an orientation reversing loop on $K$, then
$\pi_1(\mathcal F, \xi)$ is isomorphic to $\Z$.
\end{thm}

The following construction will be needed for a description of $\pi_1(\mathcal
F, \xi)$ for $\xi$ representing a
homotopically nontrivial loop on $F$ and $F\neq \R P^2, K$.

Let $F\neq \R P^2, K$ be a 
surface (not necessarily compact), and 
let $\xi$ be a curve on $F$ such that $\xi\neq 1\in\pi_1(F,
\xi(a))$. Let $f\in\pi_1(STF, \vec\xi(a))$ be the homotopy class of an
oriented fiber of the $S^1$-fibration $\pr:STF\rightarrow F$.

One can show, that there exists a unique maximal Abelian 
subgroup  $G_{\xi}<\pi_1(F, \xi(a))$ containing 
$\xi\in \pi_1(F, \xi(a))$, and that this $G_{\xi}$ is isomorphic to $\Z$ 
(see also 
Proposition~\ref{Preissman}). 
Let $g$ be its generator. Consider a curve $g_{\xi}$ direct tangent to 
$\xi$ at $\xi(a)$, which realizes $g\in\pi_1(F, \xi(a))$. 

One can show, that
$\vec\xi\in\pi_1(STF, \vec\xi(a))$ can be presented in the unique way 
as $\vec g_{\xi}^k f^l\in\pi_1(STF, \vec\xi(a))$ (see also the Proof of
Theorem~\ref{pi1nonorient}).

\begin{thm}\label{pi1nonorient}
Let $F\neq \R P^2,K$  
be a nonorientable surface (not necessarily compact) and 
let $\xi$ be a curve on $F$.

\textrm{I.} If $\xi$ represents an orientation reversing loop on $F$, 
then $\pi_1(\mathcal F,\xi)=\Z$.

\textrm{II.} If $\xi$ represents a homotopically nontrivial 
orientation preserving loop on $F$ then: 

a) $\pi_1(\mathcal F, \xi)=\Z\oplus\Z$, provided that $g_{\xi}$ is an orientation
preserving loop on $F$, or that $g_{\xi}$ is an orientation reversing loop  
and $\vec\xi=(\vec g_{\xi})^{2k}f^l$ for some nonzero $k$ and $l$.

b) $\pi_1(\mathcal F, \xi)=\pi_1(K)$, provided that $g_{\xi}$ is an
orientation reversing loop and 
$\vec\xi=(\vec g_{\xi})^{2k}$ for some nonzero $k$.

\textrm{III.} If $\xi$ represents a homotopically trivial loop, then:

a) $\pi _1(\mathcal F, \xi)$ is isomorphic to  
the subgroup of $\pi_1(STF)$ consisting of all the elements, 
which project to orientation preserving loops on $F$, provided that
$\vec\xi$ is a homotopically nontrivial loop in $STF$.
(This means, cf.~\ref{H-principle}, 
that $\xi$ is not regularly homotopic to the figure eight
curve.) 

b) $\pi_1(\mathcal F, \xi)=\pi_1(STF)$, provided that
$\vec\xi$ is a homotopically trivial loop in $STF$. (This means,
cf.~\ref{H-principle}, that $\xi$
is regularly homotopic to the figure eight curve.)
\end{thm}

The proofs of Theorems~\ref{pi1RP2},~\ref{pi1K} and~\ref{pi1nonorient} are
given in Subsections~\ref{pfpi1RP2},~\ref{pfpi1K} and~\ref{pi1arbitrary},
respectively.

\subsection{Higher homotopy groups of the space of curves.}
\begin{thm}\label{pin}
Let $F$ be a surface (not necessarily compact or orientable)
and let $\xi$ be a curve on $F$.

\textrm{I.} If $F$ is equal to $S^2$ or  $\R P^2$, then $\pi_2(\mathcal F, \xi)=\Z$ and 
$\pi_n(\mathcal F,
\xi)=\pi_n(S^2)\oplus\pi_{n+1}(S^2)$, $n\geq 3$.

\textrm{II.} If $F\neq S^2, \R P^2$, then $\pi_n(\mathcal F, \xi)=0$, $n\geq 2$. 
\end{thm}

For the Proof of Theorem~\ref{pin} see Subsection~\ref{pfpin}.

\section{Proofs}
\subsection{Some useful facts and technical Lemmas.}

\begin{lem}\label{commute} 
Let $F$ be a surface, let $STF$ be its spherical tangent bundle
and let $p\in STF$ be a point.
Let $f\in \pi_1(STF,p)$ be
the class of an oriented (in some way) 
fiber of the $S^1$-fibration $\pr: STF\rightarrow F$. 

If $\alpha\in \pi_1(STF,p)$ is a loop, which projects to an orientation
preserving loop on $F$, then 
\begin{equation}\label{commute1} 
\alpha f=f\alpha.
\end{equation}

If $\alpha\in \pi_1(STF,p)$ is a loop, which projects to an orientation
reversing loop on $F$, then 
\begin{equation}\label{commute2}
\alpha f=f^{-1}\alpha.
\end{equation}
\end{lem}

The proof of this Lemma is straightforward.

\begin{emf}[{\bf Parametric $h$-principle.}]\label{H-principle}
The parametric $h$-principle, see~\cite{Gromov} page 16, implies
that $\mathcal F$ is weak homotopy equivalent to the space $\Omega STF$ 
of free loops in $STF$. 
The corresponding mapping $h:\mathcal F\rightarrow \Omega STF$ 
sends an immersion $\xi\in \mathcal F$ 
to a loop $\vec \xi \in \Omega  STF$ by  mapping a point $y\in S^1$ 
to the point in $STF$ corresponding to the 
velocity vector of $\xi$ at $y$. 
\end{emf}

\begin{emf}\label{loopspacepi1}
{\bf Relations between the homotopy groups of the space $STF$
and of the space of free loops in $STF$.}
Let $b$ be a point in $STF$. We denote by $\Omega_b STF$ the space of all
loops in $STF$ based at $b$.

Let $\Omega STF$ be the space of all free loops in $STF$ and let $\lambda$
be a fixed element of $\Omega STF$. Fix a point $a$ on $S^1$.

Let $t:\Omega STF\rightarrow
STF$ be the mapping, which sends $\omega\in\Omega STF$ to $\omega(a)\in
STF$. One verifies that $t$ is a Serre fibration with 
the fiber isomorphic to the space of loops
based at the corresponding point.

This fibration gives rise to the following long exact sequence:

\begin{equation}\label{exact}
\cdots \stackrel{\p}{\rightarrow}
\pi_n(\Omega_{\lambda(a)}STF,\lambda)\stackrel{in_*}{\rightarrow}
\pi_n(\Omega STF, \lambda)\stackrel{t_*}{\rightarrow}
\pi_n(STF, \lambda(a))\stackrel{\p}{\rightarrow}\cdots.
\end{equation}

\end{emf}

The following statement is well known.

\begin{emf}\label{isomorphism}
Let $\lambda$ be a loop in $STF$ (not necessarily contractible)
and let $a$ be a
fixed point on $S^1$, then for any $n\geq 1:$ 
$\pi_n(\Omega_{\lambda(a)}STF,\lambda)=\pi_{n+1}(STF, \lambda(a))$.
\end{emf}

%
%
%

\begin{lem}\label{G1}(Cf. V.L.~Hansen~\cite{Hansen})
The group $\pi_1(\Omega STF, \lambda)$  is isomorphic to $Z(\lambda)$, 
the centralizer of $\lambda\in \pi_1(STF, \lambda(a))$.
\end{lem}

\begin{emf}{\bf Proof of Lemma~\ref{G1}.}
Let $t:\Omega STF\rightarrow
STF$ be the mapping described above.
A Proposition proved by V.L.~Hansen~\cite{Hansen} says that: if $X$ is a
topological space with $\pi_2(X)=0$, then
$\pi_1(\Omega X, \omega)$ is isomorphic to 
$Z(\omega)<\pi_1(X, \omega (a))$. (Here $\Omega X$
is the space of free loops in $X$ and $\omega$ is an element of $\Omega X$.)
One can verify that $\pi_2(STF)=0$ for any surface $F$. Thus, we get
that $\pi_1(\Omega STF, \lambda)$ is isomorphic to $Z(\lambda)<\pi_1(STF,
\lambda(a))$. From the proof of the Hansen's Proposition it follows that the
isomorphism is induced by $t_*$. \qed
\end{emf}

%
%
%
%
 
The following statement is an immediate consequence of Lemma~\ref{G1} and
the $h$-principle.

\begin{cor}\label{center}
Let $F$ be a surface and let $\xi$ be a curve on $F$, then $\pi_1(\mathcal F, \xi)$
is isomorphic to $Z(\vec\xi)$, the centralizer of $\vec\xi\in\pi_1(STF,
\vec\xi(a))$.
\end{cor}

\begin{lem}\label{Preissman}
Let $F\neq S^2, T^2\text{ (torus), } \R P^2, K\text{ (Klein bottle)}$ 
be a surface (not necessarily compact or orientable)  
and let $G'$ be a nontrivial commutative subgroup of $\pi_1(F)$. 
Then $G'$ is infinite cyclic 
and there exists a
unique maximal infinite cyclic group $G<\pi_1(F)$ such that $G'<G$.
\end{lem}

\begin{emf}{{\bf Proof of Lemma~\ref{Preissman}.}}
It is well known that any closed $F$, other than $S^2, T^2, \R P^2, K$,
admits a hyperbolic metric of a constant negative curvature, which is
induced from the
universal covering of it by the hyperbolic plane $H$. 
The Theorem by A.~Preissman (see~\cite{Docarmo} pp. 258-265) 
says that if $M$ is a compact Riemannian manifold with a negative curvature,
then any nontrivial Abelian subgroup $G'<\pi_1(M)$ is isomorphic to $\Z$.
Thus if $F\neq S^2, T^2, \R P^2, K$ is closed, then any nontrivial 
commutative $G'<\pi_1(F)$ is infinite cyclic.

The proof of the Preissman's Theorem given in~\cite{Docarmo} is based on
the fact, that if $\alpha,\beta\in\pi_1(M)$ are nontrivial commuting
elements, then there exists a
geodesic in $\bar M$ (the universal covering of $M$), which is mapped to itself
under the action of these elements considered as deck transformations on
$\bar M$. Moreover, these transformations restricted to the geodesic act as
translations. This implies, that if $F\neq S^2, T^2, \R P^2, K$ 
is a closed surface, then there exists a unique maximal infinite cyclic
$G<\pi_1(F)$ such that $G'<G$. This gives the proof of the Lemma for
closed $F$. 

If $F$ is not closed, then this statement is also
true, because in this case $F$ is homotopy equivalent to a bouquet of
circles.
\qed
\end{emf}

\subsection{Proof of Theorem~\ref{pi1S2}}\label{pfpi1S2}
From the exact homotopy sequence of the fibration $\pr:STS^2\rightarrow S^2$ 
we get that $\pi_1(STS^2)=\Z_2$.
Corollary~\ref{center} implies that 
$\pi_1(\mathcal F, \xi)=\Z_2$. 
\qed

\subsection{Proof of Theorem~\ref{pi1T2}}\label{pfpi1T2}
From the exact homotopy sequence of the fibration $\pr :STT^2\rightarrow T^2$
and identity~\eqref{commute1} we get that
$\pi_1(STT^2)=\Z\oplus\Z\oplus\Z$. 
Corollary~\ref{center} implies 
that $\pi_1(\mathcal F,\xi)=\pi_1(STT^2)=\Z\oplus\Z\oplus\Z$. 
\qed

\subsection{Proof of Theorem~\ref{pi1RP2}}\label{pfpi1RP2}
From the exact homotopy sequence of the fibration $\pr:ST \R P^2\rightarrow 
\R P^2$ we get that $\pi_1(ST \R P^2)=\Z_4$. 
Corollary~\ref{center} implies that $\pi_1(\mathcal F,
\xi)=\Z_4$.
\qed

\subsection{Proof of Theorem~\ref{pi1K}}\label{pfpi1K}
Corollary~\ref{center} says that 
$\pi_1(\mathcal F, \xi)=Z(\vec \xi)<\pi_1(STK, \vec\xi(a))$.

\begin{figure}[htbp]
 \begin{center}
  \epsfxsize 4cm
  \hepsffile{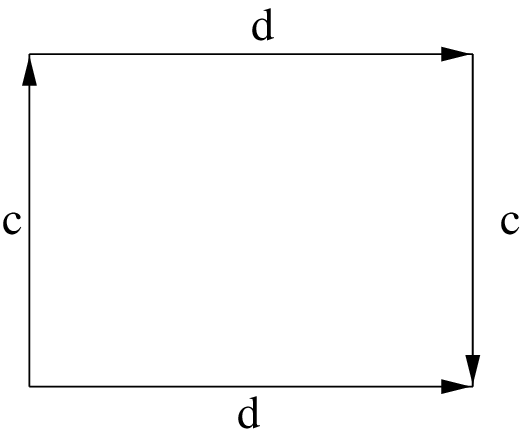}
 \end{center}
\caption{}\label{klein1.fig}
\end{figure}

Consider $K$ as a quotient of a rectangle modulo the identification on its
sides shown in Figure~\ref{klein1.fig}. We can assume that $\xi(a)$ coincides
with the image of a corner of the rectangle, and that $\xi$ and the side 
$c$ of the rectangle are
direct tangent at this point.  Let $g$ and $h$ be the
curves such that: $\vec \xi(a)=\vec g(a)=\vec h(a)$, $g=c\in \pi_1(K,
\xi (a))$ and
$h=d\in \pi_1(K, \xi (a))$. (Here $c$ and $d$ are the elements of $\pi_1(K)$
realized by the images of the sides of the rectangle used to construct $K$,
see Figure~\ref{klein1.fig}.)


One can show that:
\begin{equation}\label{present}
\pi_1(STK, \vec \xi(a))=\bigl \{\vec g, \vec h, f\big | \vec h\vec g^{\pm 1}
=\vec g^{\mp 1}\vec h, \text{ }\vec hf^{\pm
1}=f^{\mp 1}\vec h, \text{ }\vec g\ f=f\vec g\bigr \}.
\end{equation}

The second and the third relations in this presentation follow
from~\eqref{commute1} and~\eqref{commute2}.
To get the first relation one notes that
the identity $d c^{\pm 1}=c^{\mp 1} d\in\pi_1(K, \xi(a))$ implies that
$\vec h \vec g^{\pm 1}=\vec g^{\mp 1} \vec h f^k$ for some $k\in \Z$.
But $\vec h^2$ commutes with $\vec g$, since they can be lifted to $STT^2$
the fundamental group of which is Abelian. Hence $k=0$.

Relations~\eqref{present} 
on the products of $f, \vec g, \vec h$ 
imply that any element of $\pi_1(STK, \vec \xi(a))$ can be presented as 
$\vec g^k\vec h^l f ^m$ for some $k,l,m\in \Z$. Using this relations 
we calculate $Z(\vec \xi)=\pi_1(\mathcal F, \xi)$.

This group appears to be:

{\bf a)} The whole group $\pi_1(STK,\vec\xi(a))$, provided that $\vec
\xi=\vec h^{2l}$ for
some $l\in \Z$.  

{\bf b)} An isomorphic to $\Z\oplus\Z\oplus\Z$ 
subgroup of $\pi_1(STK, \vec\xi(a))$, provided that $\vec\xi=\vec g^k\vec
h^{2l}f^m$ for some $k,l,m\in\Z$ such that $k\neq 0$ or $m\neq 0$. 
This subgroup is generated by $\{\vec g, \vec h^2,f\}$.

{\bf c)} An isomorphic to $\Z$ subgroup of $\pi_1(STK, \vec\xi(a))$, 
provided that $\vec \xi=\vec g^k\vec h^{2l+1} f^m$ for some $k,l,m\in\Z$.
This subgroup is generated by $\alpha_{\xi}=\vec g^k\vec h f^m$.
(Note that $\alpha_{\xi}^2=\vec h^2$, and $\vec \xi=\alpha_{\xi}^{2l+1}$.)

Now the statement of the Theorem is a direct consequence of
relations~\eqref{present}. 
\qed

\subsection{Proof of Theorem~\ref{pi1orient} and
Theorem~\ref{pi1nonorient}.}\label{pi1arbitrary}

We are going to prove that the statement of Theorem~\ref{pi1nonorient}  is
true for any orientable surface $F\neq S^2, T^2$ and any
nonorientable $F\neq \R P^2, K$. (We will see that
$\vec\xi\in\pi_1(STF, \vec\xi(a))$ can be presented in the unique way 
as $\vec g_{\xi}^k f^l\in\pi_1(STF, \vec\xi(a))$ for any $F\neq S^2, \R P^2,
T^2, K$.)

Clearly this gives a proof of Theorem~\ref{pi1nonorient}.
Theorem~\ref{pi1orient} is also an immediate consequence of this fact.

\begin{emf}\label{xineq0}
{\bf Proof of Theorem~\ref{pi1orient} and
Theorem~\ref{pi1nonorient} in the case of $\xi\neq 1\in\pi_1(F, \xi(a))$.}
Consider a subgroup $G'$ of $\pi_1(F, \xi(a))$ generated by $\xi$. It is an
infinite cyclic group (see~\ref{Preissman}). 
There is a unique (see~\ref{Preissman}) maximal infinite cyclic 
group $G<\pi_1(F, \xi(a))$ such that $G'<G$.
Let $g$ be the generator of $G$. Let $g_{\xi}$ be a curve direct tangent to 
$\xi$ at $\xi(a)$ representing this $g$.

Corollary~\ref{center} says that $\pi_1(\mathcal F, \xi)$ is isomorphic to 
$Z(\vec \xi)$.
Take $\alpha\in Z(\vec\xi)$.
Since $\vec\xi$ and $\alpha$ commute in
$\pi_1(STF,\vec\xi(a))$, we get that their images under the projection
$\pr_*:\pi_1(STF,\vec\xi(a))\rightarrow\pi_1(F, \xi(a))$ commute in
$\pi_1(F,\xi(a))$. 
Lemma~\ref{Preissman} implies that these projections are in the subgroup
$G$.

The kernel of the homomorphism $\pr_*$ is generated by $f$, 
the homotopy class of an oriented fiber of the $S^1$-fibration
$\pr:STF\rightarrow F$. This fact and identities~\eqref{commute1} and 
\eqref{commute2} show that there exist unique $k,l,m,n\in \Z$ such that 
$\vec\xi=\vec g_{\xi}^kf^l$ and $\alpha=\vec g_{\xi}^mf^n$.

Using identities~\eqref{commute1} and~\eqref{commute2} we can check for which
values of $k,l,m,n$ the elements $\alpha$ and $\vec\xi$ commute. 
This allows us to calculate 
$Z(\vec\xi)$. It turns out to be:

a) A group isomorphic to $\Z\oplus\Z$ generated by 
$\{\vec g_{\xi}, f\}$, provided that $g_{\xi}$ is an orientation preserving loop
on $F$. 

b) A group isomorphic to $\Z$ generated by $\vec g_{\xi}f^l$, 
provided that $g_{\xi}$ is an orientation reversing loop, and $k$ is odd.
(This means that $\xi$ represents an orientation reversing loop on $F$.)
Note also, that in this case $(\vec g_{\xi}f^l)^2=\vec g_{\xi}^2$.

c) A group isomorphic to $\Z\oplus\Z$ generated by $\{\vec g_{\xi}^2, f\}$,
provided that $g_{\xi}$ is an orientation reversing loop on $F$, $k\neq 0$ is
even, and $l\neq 0$.

d) A group isomorphic to $\pi_1(K)$ generated by $\{\vec g_{\xi}, f\}$,
provided that $g_{\xi}$ is an orientation reversing loop on $F$, $k\neq 0$ is
even, and $l=0$.

(Note that if $k=0$, then $\xi=1\in\pi_1(STF, \xi(a))$, which contradicts to
our assumption.)

This finishes the proof of the two theorems for this case.
\end{emf}

\begin{emf}\label{xieq0}
{\bf Proof of Theorem~\ref{pi1orient} and
Theorem~\ref{pi1nonorient} in the case of $\xi=1\in\pi_1(F,
\xi(a))$.}
From the exact homotopy sequence of the $S^1$-fibration 
$\pr:STF\rightarrow F$ 
we get that $\ker \pr_*$ is generated by $f$, the homotopy class of the
fiber.
Since $\xi=1\in\pi_1(F, \xi(a))$ we get that there exists 
a $k\in\Z$ such that $\vec\xi=f^k$.
Lemma~\ref{split} says that 
$\pi_1(\mathcal F, \xi)$ is
isomorphic to $Z(\vec\xi)=Z(f^k)<\pi_1(STF, \vec \xi(a))$. 

For $k\neq 0$ identities~\eqref{commute1} and~\eqref{commute2} imply that 
$Z(f^k)$ coincides with the set of elements of $\pi_1(STF, \vec\xi(a))$,
which project to orientation preserving loops on $F$.
This finishes the proof of the Theorem for $\xi=1\in\pi_1(F, \xi(a))$ and 
$\vec\xi\neq 1\in \pi_1(STF, \vec\xi(a))$.

If $k=0$, then $\vec\xi=1\in \pi_1(STF, \vec\xi(a))$. Thus
$Z(\vec\xi)=\pi_1(STF, \vec\xi(a))$. Hence, in this case $\pi_1(\mathcal F,
\xi)=\pi_1(STF, \vec\xi(a))$.
\end{emf}

\subsection{Proof of Theorem~\ref{pin}.}\label{pfpin}
The proof of this Theorem is based on the following exact sequence, which
was introduced in section~\ref{loopspacepi1}.

\begin{equation}\label{temp}
\cdots \stackrel{\p}{\rightarrow}
\pi_n(\Omega_{\lambda(a)}STF,\lambda)\stackrel{in_*}{\rightarrow}
\pi_n(\Omega STF, \lambda)\stackrel{t_*}{\rightarrow}
\pi_n(STF, \lambda(a))\stackrel{\p}{\rightarrow}\cdots.
\end{equation}

\begin{lem}\label{split}
If $F$ is equal to $S^2$ or $\R P^2$ and $n\geq 2$, then  
$$\pi_n(\Omega STF, \lambda)=\pi_n(\Omega_{\lambda(a)}STF, \lambda)\oplus
\pi_n(STF, \lambda(a)).$$
\end{lem}

\begin{emf}{\bf Proof of Lemma~\ref{split}}
Fix $n>1$. We construct a homomorphism 
$g:\pi_n(STF, \lambda(a))\rightarrow\pi_n(\Omega STF, \lambda)$ such that 
$t_*\circ g =\id_{\pi_n(STF, \lambda(a))}$. After this the exactness of the
sequence~\eqref{temp} and the fact that higher homotopy groups are
Abelian imply the statement of the Lemma.

We describe this construction for $F=\R P^2$. The construction of $g$ for
$F=S^2$ can be easily deduced from this one.

From the exact homotopy sequence of the covering $ST S^2\rightarrow ST\R
P^2$ we get that $\pi_n(ST \R P^2)$, $n\geq 2$, 
is canonically isomorphic to 
$\pi_n(ST S^2)$. 

Take $s:S^n\rightarrow ST \R P^2$, which represents a given element 
of $\pi_n(ST \R P^2, \lambda(a))$. 
Let $s':S^n\rightarrow STS^2$ be the mapping which
is a lifting of $s$ under the covering $STS^2\rightarrow ST\R P^2$.  
Fix an orientation of $S^2$. 
Then for every $x\in S^n$ the orientation of a small neighborhood of 
$\pr s'(x)\in S^2$ induces an orientation of a small
neighborhood of $\pr s(x)\in \R P^2$. 

There is a unique isometric autodiffeomorphism 
$I_x$ of $\R P^2$ such that:

a) It maps $\pr s(*)$ to $\pr s(x)$.

b) The differential of it sends $s(*)$ to $s(x)$.

c) The above described local orientation at $\pr s(x)$ coincides with the
one induced by the differential of $I_x$ from the local orientation at 
$\pr s(*)$. 

Let $\bar s:S^n\rightarrow \Omega ST \R P^2$ be the mapping
which sends $x\in S^n$ to $I_x(\lambda)$ (the translation of $\lambda$ by
$I_x$).

Set the value of $g$ on the element of 
$\pi_n(ST\R P^2, \lambda(a))$ represented by $s$ to be the element of 
$\pi_n(\Omega ST\R P^2, \lambda)$ represented by $\bar s$. 
A straightforward verification
shows that this $g$ is the desired 
homomorphism from $\pi_n(ST \R P^2, \lambda(a))$ 
to $\pi_n(\Omega ST\R P^2, \lambda)$.

This finishes the proof of Lemma~\ref{split}.
\qed
\end{emf}

\begin{emf}
One verifies that 
$\pi_2(STF)=0$ and $\pi_n(STF)=\pi_n(S^2)$, $n\geq 3$, for $F$ equal to $S^2$
or $\R P^2$. Now Lemma~\ref{split}, statement~\ref{isomorphism} and the
weak homotopy equivalence given by the 
$h$-principle (see~\ref{H-principle}) imply the
first statement of the Theorem. (Note that 
$\pi_3(S^2)=\Z$.)  

Statement~\ref{isomorphism} says that 
$\pi_n(\Omega_{\lambda(a)}STF, \lambda)=\pi_{n+1}(STF, \lambda(a))$.
One verifies that $\pi_n(STF)=0$, $n\geq2$ for $F\neq S^2, \R P^2$.
The exactness of sequence~\eqref{temp} implies, that $\pi_n(\Omega STF,
\lambda)=0$, $n\geq 2$. Using the weak homotopy equivalence 
given by the $h$-principle we get the second
statement of the Theorem.

This finishes the proof of Theorem~\ref{pin}.
\qed 
\end{emf}

\bibliographystyle{amsplain}

\end{document}